\documentclass[11pt]{amsart}

\usepackage[paperwidth=20cm,paperheight=29.7cm,left=2cm,top=2cm]{geometry}

\usepackage[alphabetic]{amsrefs}
\usepackage{graphicx}
\graphicspath{ {images/} }
\usepackage[T1]{fontenc}
\usepackage{txfonts}
\usepackage{times}
\usepackage{amssymb,verbatim,amscd,amsfonts,amsmath,amsthm,enumerate}
\usepackage[all]{xy}
\usepackage{subfig}
\usepackage{tikz}
\usetikzlibrary{shapes,arrows,shadows}
\usetikzlibrary{decorations.markings}
\usepackage{color}
\usepackage[normalem]{ulem}

\usepackage{hyperref}
\usepackage{wrapfig}
\usetikzlibrary{arrows}

\long\def\symbolfootnote[#1]#2{\begingroup%
\def\thefootnote{\fnsymbol{footnote}}\footnote[#1]{#2}\endgroup}

\newtheorem{theorem}{Theorem}[section]

\theoremstyle{definition}

\newtheorem{remark}[theorem]{Remark}

\newtheorem{definition}[theorem]{Definition}

\newcommand{\N}{\mathbb{N}}

\begin{document}
\title{On infinite regular and chiral maps}

\author{John A. Arredondo, Camilo Ram\'irez and Ferr\'an Valdez}

\address{
Universidad Konrad Lorenz, Bogot\'a Colombia, C.P. 110231, Centro de Ciencias Matem\'aticas, UNAM, Campus Morelia, C.P. 58190, Morelia, Michoac\'an,  M\'exico.
}
\email{alexander.arredondo@konradlorenz.edu.co, camilomaluendas@gmail.com,\newline ferran@matmor.unam.mx}

\begin{abstract}
We prove that infinite regular and chiral maps take place on surfaces with at most one end. Moreover, we prove that an infinite regular or chiral map on an orientable surface with genus can only be realized on the Loch Ness monster, that is, the topological surface of infinite genus with one end.
\end{abstract}

\maketitle


\section{Introduction}
	\label{SEC:Introduction}
This paper is motivated by the following problem, posed by D. Pellicer in 	\cite{Pellicer12}: \emph{determine which noncompact surfaces without boundary admit embeddings of chiral maps}. The present work gives a complete answer to this question and generalizes previous results on minimal regular covers of the Archimedean tessellations, see \cite{CPRRV15}. More precisely, we prove the following  four theorems:

\begin{theorem}\label{TH:EndsCayley}
 Let $\mathcal{M}$ be a regular or chiral map on a surface $S$, and let $Aut(\mathcal{M})$ be the automorphism group of the map $\mathcal{M}$. Then the spaces $Ends(S)$ and $Ends(Aut(\mathcal{M}))$ are in bijection.
\end{theorem}


\begin{theorem}
	\label{TH:endsAutGroups}
Let $\mathcal{M}$ be a regular or chiral map. Then the cardinal of $Ends(Aut(\mathcal{M}))$ is either zero or one.
\end{theorem}

Here $Ends(S)$ denotes the space of ends of the surface $S$ and $Ends(G)$ the space of ends of the Cayley graph of $G$ with respect some generator subset $H\subset G$. For a precise definition of these spaces see \S\ref{SEC:Preliminaries}. The preceding theorems tell us that there is a considerable topological restriction for a surface to support an infinite chiral or regular map, namely, it necessarily has to have one-ended. As we will see in \S \ref{SUBSEC:ProofCorollaries}, we can describe the whole topological zoology of infinite surfaces supporting orientable regular and chiral maps:



\begin{theorem}
	\label{Coro:RegularMapsMonster}
Let $\mathcal{M}$ be a regular map on a noncompact and orientable surface $S$. Then $S$ is homeomorphic to either the plane or the Loch Ness monster.
\end{theorem}

\begin{theorem}
	\label{Coro:ChiralMapsMonster}
Let $\mathcal{M}$ be a chiral map on a noncompact surface $S$. Then $S$ is homeomorphic to the Loch Ness monster.
\end{theorem}

The Loch Ness monster is the only orientable topological surface with infinite genus and only one end. As a mathematical object, this surface appears naturally in many contexts, see \textit{e.g.}, \cite{Cox1}, \cite{Valdez1}, \cite{Valdez2}, and \cite{Ghys}.\\

\textbf{Acknowledgments}. The authors thank D. Pellicer for helpful conversations and H. Wilton for answering one of our questions via Mathoverflow \cite{Valdez3}. The second author was partially supported by CONACYT, CCM-UNAM, and IFM-UMSNH. The last author was generously supported by CONACYT CB-2009-01
127991 and PAPIIT projects IN100115, IN103411 \& IB100212 during the realization of this project.


\section{Preliminaries}
	\label{SEC:Preliminaries}

\textbf{Maps}. We begin this section by discussing some general aspects of maps that will be needed for the proofs of the main theorems. Our text is not self-contained, hence we refer the reader to \cite{Pellicer12}, \cite{JoSi} and references within for details.

A \textbf{map} $\mathcal{M}$ on a surface $S$ is a 2-cell embedding $i:\Gamma\hookrightarrow S$ of  \textbf{locally finite simple} graph\footnote{We think of $\Gamma$ as the geometric realization of an abstract graph.}  $\Gamma$ into the surface $S$. In other words, only finitely many edges are incident on each vertex of $\Gamma$, the ends of each edge are always different vertices and the function $i$ is a topological embedding such that each connected component of $S\setminus i(\Gamma)$ is homeomorphic to a disc. We denote such a triple by $\mathcal{M}:=\mathcal{M}(\Gamma,i,S)$.

Every $f\in Homeo(S)$ for which there exists an isomorphism $g_f:\Gamma\to\Gamma$ such that $i\circ g_f = f\circ i$, in other words such the following diagram is commutative, is called a \textbf{preautomorphism} of the map $\mathcal{M}(\Gamma,i,S)$.

\[
\xymatrix{S\ar[rr]^f & & S\\ \Gamma\ar@{^(->}[u]^{i}\ar[rr]_{g_f}& & \Gamma\ar@{^(->}[u]_{i}}
\]

The set of preautomorphisms of a map $\mathcal{M}$ has a natural group structure and we denote it by $\widetilde{Aut}(\mathcal{M})$. It is not difficult to see that for each $f\in \widetilde{Aut}(\mathcal{M})$ the map $g_f$ is unique, hence we have a group morphism $\varphi:\widetilde{Aut}(\mathcal{M})\to Isom(\Gamma)$, where the codomain is the group of isomorphism of $\Gamma$. An \textbf{automorphism of the map} $\mathcal{M}$ is an element of the group $Aut(\mathcal{M}):=\widetilde{Aut}(\mathcal{M})/{\rm Ker\varphi}$.
 We shall abuse notation and we shall write the coset $[f]\in Aut(\mathcal{M})$ as $f$.

\begin{figure}[h]
    \centering
    \includegraphics[scale=0.44]{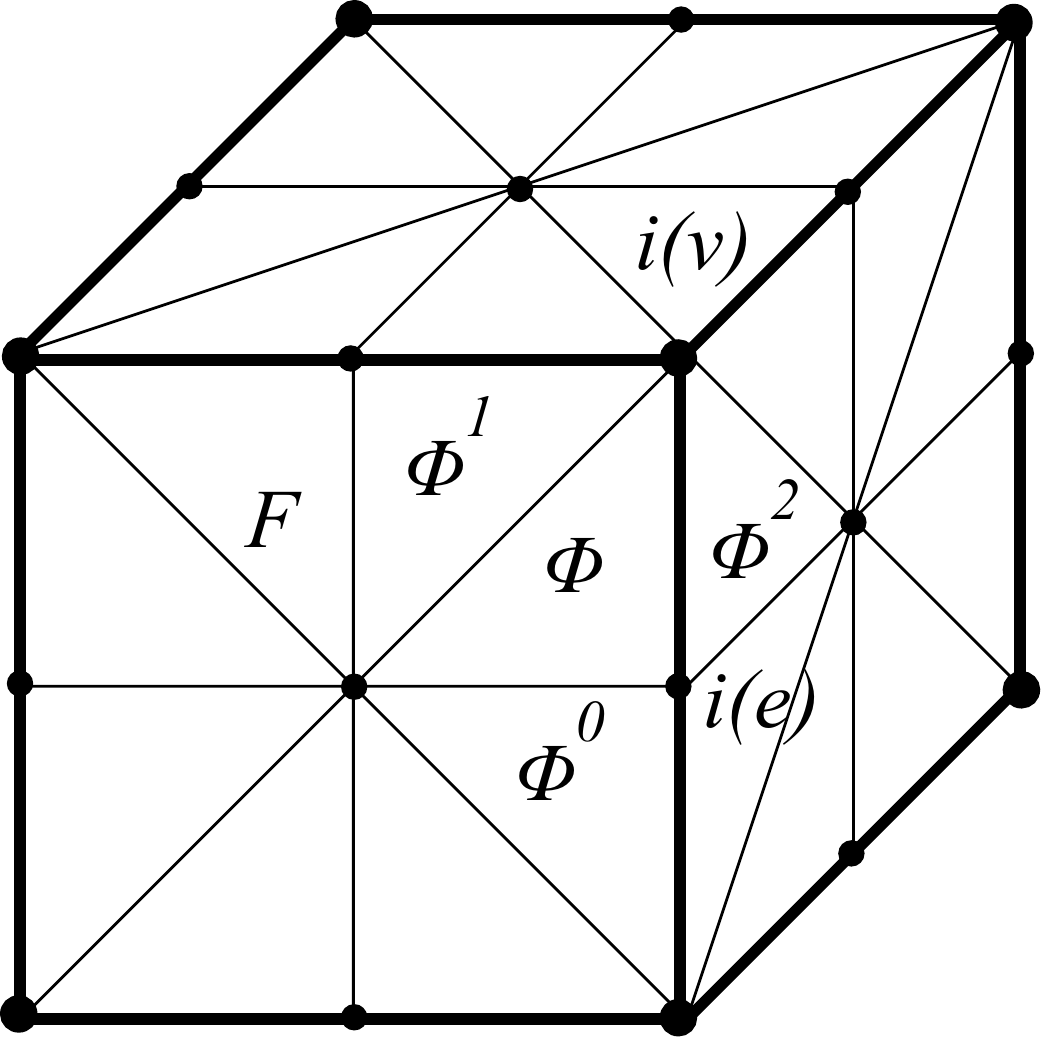}\\
    \caption{A map in a cube divided into flags.}
    \label{FIGURE1}
\end{figure}

A \textbf{flag} $\Phi$ of a map $\mathcal{M}(\Gamma,i,S)$ is a triangle on $S$ whose vertices are a \textbf{vertex} $i(v), v\in V(\Gamma)$, the midpoint of an \textbf{edge} $i(e)$ containing $i(v)$ with $e\in E(\Gamma)$, and an interior point of a \textbf{face} $F\subset S\setminus i(\Gamma)$ whose boundary contains $i(e)$. We may assume that for all flags containing the closure in $S$ of the face $F$, the same interior point of $F$ is chosen to be a vertex of the corresponding triangle. This way, every map induces a triangulation of $S$ given by its flags. Combinatorially, every flag can be identified with an ordered incident triple of vertex, edge and face of the map $\mathcal{M}$. Henceforth, we shall abuse notation and we shall understand flags either as triangles or as ordered triples. Given a flag $\Phi$ of the map $\mathcal{M}$, there is a unique adjacent flag $\Phi^0$ (resp. $\Phi^1$ and $\Phi^2$) of the map $\mathcal{M}$ that differs from $\Phi$ precisely on the vertex (resp. on the edge and on the face). The flag $\Phi^j$ is called the $j$-\textbf{adjacent flag of} $\Phi$. In figure \ref{FIGURE1}, we show an example of the cube with some flags marked with their name.
We denote the set of flags of a given map $\mathcal{M}$ by $\mathcal{F}:=\mathcal{F}(\Gamma,i,S)$. The group $Aut(\mathcal{M})$ acts on the set of flags $\mathcal{F}(\Gamma,i,S)$ and for every pair $\Phi_1, \Phi_2$ of flags of $\mathcal{M}$, there exist at most one automorphism of $\mathcal{M}$ mapping $\Phi_1$ to $\Phi_2$. In other words, every elements of $Aut(\mathcal{M})$ is completely determined by the image of a given flag.

\begin{definition}[Regular and Chiral maps]
	\label{DEF:REGANDCHIRAL}
A map $\mathcal{M}$ is called \textbf{regular}, respectively \textbf{chiral}, if the action of $Aut(M)$ on $\mathcal{F}$ induces one orbits in flags, respectively the action of $Aut(M)$ on $\mathcal{F}$ induces two orbits in flags with the property that adjacent flags belong to different orbits.
\end{definition}

The graph $\Gamma$  of a regular or chiral map $\mathcal{M}(\Gamma,i,S)$ is always regular, that is, every vertex  has the same degree $q\in\mathbb{N}$. Moreover, such maps also satisfy that every boundary of each connected component of $S\setminus i(\Gamma)$ is formed by a cycle in $i(\Gamma)$ of fixed length $p$. The pair $\{p,q\}$ is called the \textbf{Schl\"afli type} of the map $\mathcal{M}$. When $\mathcal{M}$ is a regular map the group $Aut(\mathcal{M})$ is generated by three involutions $\rho_{0},\rho_{1},$ and $\rho_{2}$, where $\rho_{j}$ is the unique automorphism of $\mathcal{M}$ sending a fixed \textbf{base flag} $\Phi$ to its $j$-adjacent flag $\Phi^{j}$. Moreover the generating set $\{\rho_{0},\rho_{1},\rho_{2}\}$ satisfies the relations:
\begin{equation}
 	\label{E:RELATIONS AUTOMGROUP}
\rho_{1}^{2}=\rho_{1}^{2}=\rho_{1}^{2}=(\rho_{0}\rho_{2})^{2}=(\rho_{0}\rho_{1})^{p}=(\rho_{1}\rho_{2})^{q}=Id,
\end{equation}
and probably other more. On the other hand, when $\mathcal{M}$ is a chiral map the group $Aut(\mathcal{M})$ is generated by two elements $\{R,S\}$ that satisfy the relations:
\begin{equation}
R^{p}=S^{q}=Id,
\end{equation}
and probably other more. Roughly speaking, the element $R$ acts on the set of flags of the map $\mathcal{M}$ as a rotation of order $p$ centered in the middle of a given face of the map $\mathcal{M}$. Analogously, if $\Phi_1$ and $\Phi_2$ are different flags of the chiral map $\mathcal{M}$ but they belong to same orbit, then there is a unique automorphism of $\mathcal{M}$ mapping $\Phi_{1}$ to $\Phi_{2}$.

Finally let us mention that no chiral map lies on a non-orientable surface whereas, on the contrary, there are many non-orientable surfaces with regular maps.

\vspace{2mm}
\textbf{The space of ends}. We start by discussing the space of ends of a topological space $X$. We perform this in full generality, though we shall only use it when $X$ is a surface or a graph.

Let $X$ be a locally compact, locally connected, connected, and Hausdorff space. There are several equivalent ways to define the space of ends of $X$, the following was introduced by Freudenthal.

\begin{definition}\cite{Freu31}
Let $U_{1}\supseteq U_{2}\supseteq\ldots$ be an infinite nested sequence of non-empty connected open subsets of $X$ such that for each $n\in\N$ the boundary $\partial U$ of $U_{n}$ is compact, $\bigcap\limits_{n\in\N}\overline{U_{n}}=\emptyset$ and for any compact  $K\subset X$ there exist $m\in\mathbb{N}$ such that $K\cap U_m =\emptyset$. We denote the sequence $U_{1}\supseteq U_{2}\supseteq\ldots$ as $(U_n)_{n\in\mathbb{N}}$. Two such sequences $(U_n)_{n\in\mathbb{N}}$ and $(U'_{n})_{n\in\mathbb{N}}$ are said to be equivalent if for every $l\in\N$ there exist $k$ such that $U_{l}\supseteq U'_{k}$ and viceversa, that is, for every $n\in\N$ there exist $m$ such that $U'_{n}\supseteq U_{m}$. The corresponding equivalence classes are also called \emph{topological ends} of $X$ and we will denote it by ${\rm Ends}(S)$.
\end{definition}

For every non-empty open subset $U$ of $X$ such that its boundary $\partial U$ is compact we define
\begin{equation}\label{E:BaseTopEnds}
U^{*}:=\{[U_{n}]_{n\in\mathbb{N}}\in{\rm Ends}(X)\hspace{1mm}|\hspace{1mm}U_{j}\subset U\hspace{1mm}\text{for some }j\in\mathbb{N}\}.
\end{equation}
The collection formed by all sets of the form $U\cup U^{*}$, with $U$ open with compact boundary of $X$, forms a base for the topology of $X':=X\cup{\rm Ends}(X)$.

\begin{theorem}\cite{Ray60}
	\label{TH:TOPOLENDS}
Let $X':=X\cup{\rm Ends}(X)$ be the topological space defined above. Then,
\begin{enumerate}
\item The space $X'$ is Hausdorff, connected and locally connected.
\item The space ${\rm Ends}(X)$ is closed and has no interior points in $X'$.
\item The space ${\rm Ends}(X)$ is totally disconnected in $X'$.
\item The space $X'$ is compact.
\item If $V$ is any open connected set in $X'$, then $V\setminus {\rm Ends}(X)$ is connected.
\end{enumerate}
\end{theorem}

\textbf{Ends of surfaces}. When $X$ is a surface $S,$ the space ${\rm Ends}(S)$ carries extra information, namely, those ends that carry \textbf{infinite} genus. This data, together with the genus and the orientability class, determine the topology of $S$. We discuss the details of this fact in the following paragraphs. Given that in this article we only deal with orientable surfaces, from now on we dismiss the non-orientable case.

The topological classification of finite type compact topological surfaces is well known. Indeed, the topological invariant in this case is just a pair of non negative integers $(g,b)$ denoting genus and number of boundary components respectively, plus the orientability class of the surface.
A surface is said to be \textbf{planar} if all of its compact subsurfaces are of genus zero. An end $[U_n]_{n\in\mathbb{N}}$ is called \textbf{planar} if there exists
an $l\in\N$ such that $U_l$ is planar.

\begin{definition}
The \textbf{genus} of a surface $S$ is the maximum of the genera of its compact subsurfaces (with boundaries).
\end{definition}
Remark that if a surface $S$ has \textbf{infinite genus} there exists no finite set $\rm \mathcal{C}$ of mutually non-intersecting simple closed curves with the property that $S\setminus \mathcal{C}$ is  \textbf{connected and planar}. We define ${\rm Ends_{\infty}}(S)\subset{\rm Ends}(S)$ as the set of all ends of $S$ which are not planar. It follows from the definitions that  ${\rm Ends_{\infty}}(S)$ forms a closed subspace of  ${\rm Ends}(S)$. 

\begin{theorem} (Classification of orientable surfaces, \cite[Chapter 5]{Kerek})
Two orientable surfaces $S$ and $S'$  are homeomorphic if and only if they have the same genus $g\in\mathbb{N}\cup\{0,\infty\}$, and both ${\rm Ends_{\infty}}(S)\subset {\rm Ends}(S)$ and ${\rm Ends_{\infty}}(S')\subset {\rm Ends}(S')$ are homeomorphic as nested topological spaces.
\end{theorem}

\begin{definition}
\cite{PhiSul} The surface with one end and infinite genus is called the \textit{Loch Ness monster} (see Figure 2).
\begin{figure}[h]\label{FIGURE2}
    \centering
    \includegraphics[scale=0.7]{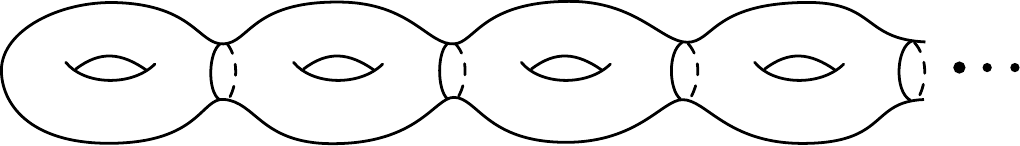}\\
    \caption{Loch Ness monster.}
\end{figure}
\end{definition}

\section{Proof of main results}
In this section we present the proofs of our main results.

\subsection{Proof theorem \ref{TH:EndsCayley}} We divide the proof in two cases.\\
\\
\textbf{Case 1: regular maps}. Let $\mathcal{M}=\mathcal{M}(\Gamma,i,S)$ be a regular map and $G$ the Cayley graph of $Aut(\mathcal{M})$ with respect to the generator set $\{\rho_0,\rho_1,\rho_2\}$. The idea of the proof is to define an embedding $i:G\hookrightarrow S$ and then prove that it induces a bijection $i_{\ast}:Ends(G)\to Ends(S)$.

\begin{figure}[h]
    \centering
    \includegraphics[scale=0.45]{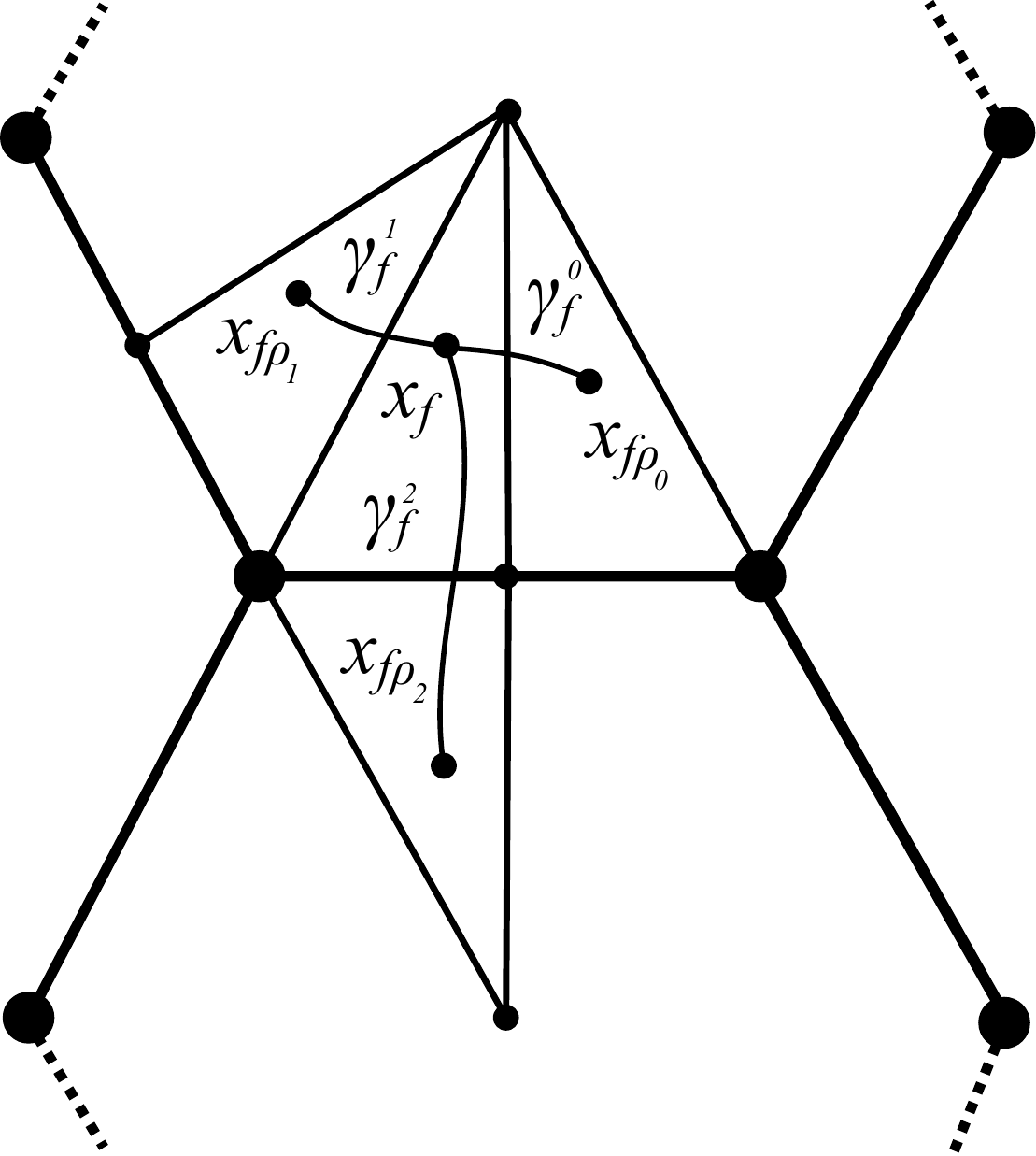}\\
    \caption{Embedding $i:G\hookrightarrow S$.}
    \label{FIGURE3}
\end{figure}

\vspace{2mm}
\emph{The embedding}. Fix a base flag $\Phi_{Id}$ of the map $\mathcal{M}$. Given that a flag $\Phi$ is a triangle on $S$ and the action of $Aut(\mathcal{M})$ on $\mathcal{F}$ is transitive we have that:
\begin{enumerate}
\item For each $\Phi\in\mathcal{F}\setminus \{\Phi_{Id}\}$ there is a unique element $f\in Aut(\mathcal{M})$ such that $f(\Phi_{Id})=\Phi$. We relabel the flag $\Phi$ as $\Phi_f$ and pick for each $f\in Aut(\mathcal{M})$ a point $x_f$ belongs to the interior of the flag $\Phi_f$.
\item For each flag $\Phi_f$ and its respective j-adjacent flags $\Phi_{f\rho_j}$, $j\in \{0,1,2\}$, we can choose three simple paths  $\{\gamma_f^j\}_{j=0}^2$ in $S$ such that $\gamma_f^j$ has end points $x_f$ and $x_{f\rho_j}$, $\cap_{j=0}^2\gamma_f^j=x_f$ and $\gamma_f^j=\gamma_{f\rho_j}^j$ for each $j\in \{0,1,2\}$.
\end{enumerate}
In figure \ref{FIGURE3} we illustrate the definition of $i:G\hookrightarrow S$.
As it can be seen, on vertices $i(f)=x_f$ and on edges $i((f,f\rho_j))=\gamma_f^j$, for each $j\in \{0,1,2\}$. It is clear that our choice of $i$ is a topological embedding.

\vspace{2mm}
Let $U$ be a connected open  subset of $G$ with compact boundary, and  without loss of generality, such that $\partial U \subset V(G)$.
Define:
\begin{equation}
\widetilde{U}:= Int \left(\bigcup_{f\in U\cap V(G)} \Phi_f\right)\subset S.
\end{equation}
In what follows we show that:
\begin{align}
	\label{E:HomeoEnds}
i_{\ast}: Ends(G)&\to Ends(S)\\
[U_n]_{n\in \N}&\to [\widetilde{U}_n]_{n\in \N},\nonumber
\end{align}
is a (well defined) bijection.

\begin{remark}\label{ob_2}
\begin{upshape}
Let $U_1, U_2\subset G$ be connected open subsets with compact boundary such that  $\partial U_1, \partial U_2 \subset V(G)$. Then
\begin{enumerate}

\item The inclusion $U_1\supset U_2$ implies $\widetilde{U}_1\supset \widetilde{U}_2$.
\item The intersection $U_1\cap U_2=\emptyset$ implies $\widetilde{U}_1\cap \widetilde{U}_2=\emptyset$.
\end{enumerate}
\end{upshape}
\end{remark}
From this implications we deduce that if $(U_n)_{n\in\N}$ is a representative of $[U_n]_{n\in\N}$, then $(\widetilde{U}_n)_{n\in\N}$ is a nested sequence. Now for every compact $K\subset S$ we define the compact
\[
\widetilde{K}:=\{f\in V(G): \Phi_f \cap  K\neq\emptyset\}\subset G.
\]
Since $[U_n]_{n\in \N}$ is an end of $G$ there exist $m\in\N$ such that $U_m\cap \widetilde{K}=\emptyset$ and therefore $\widetilde{U}_m\cap K=\emptyset$. This proves that $[\tilde{U}_n]_{n\in\mathbb{N}}$ is an end of $S$.The second implication of remark \ref{ob_2} gives us the well-definiteness of  $i_{\ast}$. If the ends $[U_n]_{n\in \N}$ and $[V_n]_{n\in \N}$ of $G$ are not equivalent, then there exist $n\in \N$ such that $U_n \cap V_n=\emptyset$. By Remark \ref{ob_2} we have that $\widetilde{U}_n \cap \widetilde{V}_n=\emptyset $, \textit{i.e.},  the ends $[\widetilde{U}_n]_{n\in \N}$ and $[\widetilde{V}_n]_{n\in \N}$ of $S$ are no equivalent. In other words, the map $i_{\ast}$ is injective. Moreover, if  $[\widetilde{U}_{n}]_{n\in\mathbb{N}}$ is an element in $Ends(S)$, then for all $n\in\mathbb{N}$ $U_n$ is a component connected of $i^{-1}(G\cap \widetilde{U})$ such that $U_n\supset U_{n+1}$. This proves that $i_{\ast}$ is surjective. Hence $i_{\ast}[U_{n}]_{n\in\mathbb{N}}=[\widetilde{U}_{n}]_{n\in\mathbb{N}}$. This implies that $i_{\ast}$ is a bijection.

\vspace{2mm}
\textbf{Case 2: chiral maps}. Let $\mathcal{M}=\mathcal{M}(\Gamma,i,S)$ be a chiral map and $G$ the Cayley graph of $Aut(\mathcal{M})$ with respect to the generator set of rotations $\{R,S\}$. The idea of the proof is the same as in the regular case: first we define an embedding $i:G\hookrightarrow S$ and then we prove that it induces a bijection $i_\ast:Ends(G)\to Ends(S)$. The only main difference is how the embedding $i$ is defined, and this we explain in the next paragraphs.

\begin{figure}[h]
    \centering
    \includegraphics[scale=0.45]{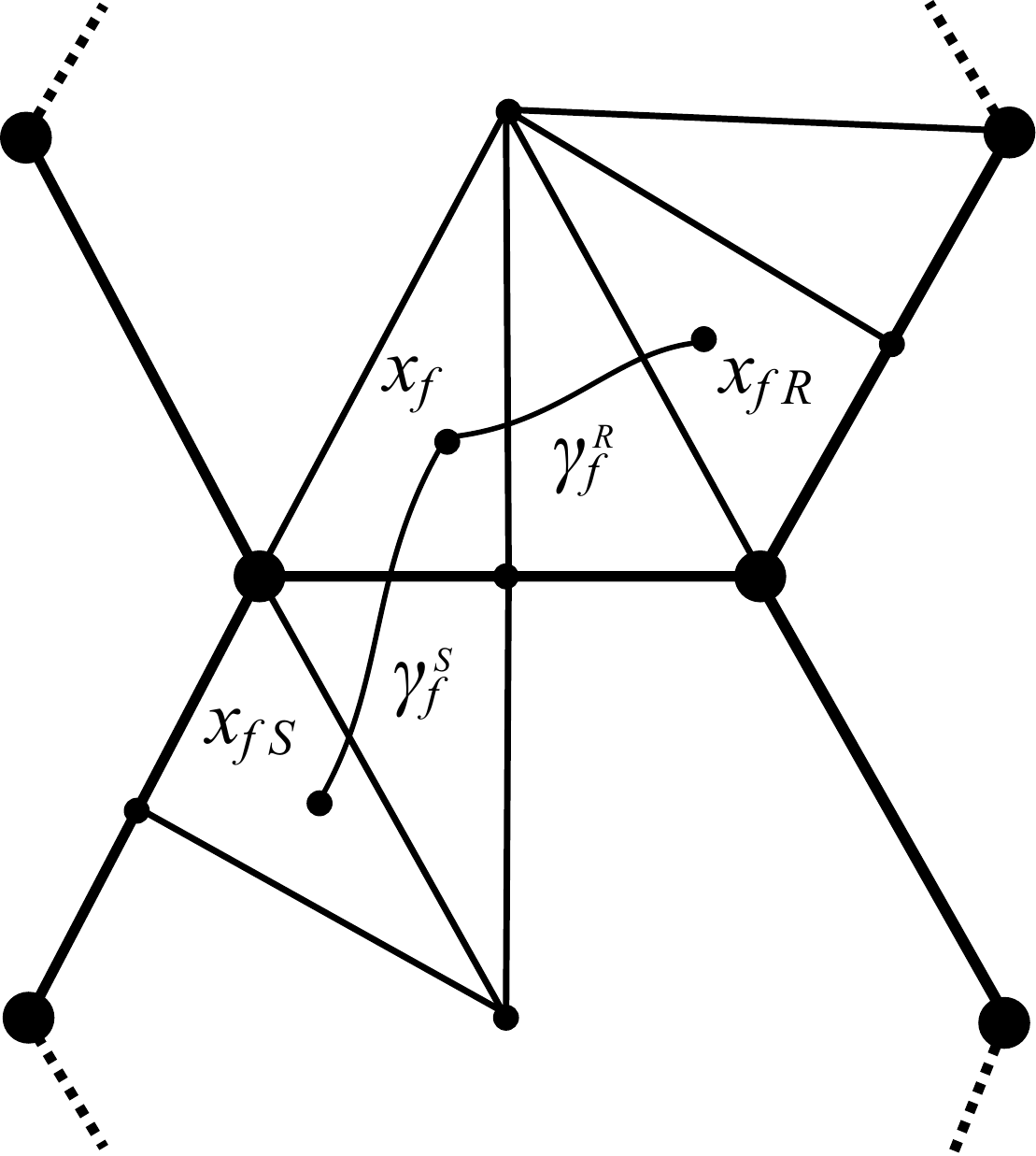}\\
    \caption{Embedding $i:G\hookrightarrow S$.}
    \label{FIGURE4}
\end{figure}

\vspace{2mm}
\emph{The embedding}. Let $[\Phi]$ and $[\Psi]$ denote the two orbits of the $Aut(\mathcal{M})$-action on $\mathcal{F}$ and pick a base flag $\Phi_{Id}\in [\Phi]$. Let $\Psi_{Id}\in [\Psi]$ be the unique adjacent flag to $\Phi_{Id}$ that differs on a vertex. Given that the restriction of the $Aut(M)$-action on $\mathcal{F}$ to $[\Phi]$ or $[\Psi]$ is transitive we have that:
\begin{enumerate}
\item For each $\Phi\in[\Phi]\setminus \{\Phi_{Id}\}$ there is a unique element $f \in Aut(\mathcal{M})$ such that $f(\Phi_{Id})=\Phi$. We relabel $\Phi$ as $\Phi_f$ and pick for each $f\in Aut(\mathcal{M})$ an interior point $x_f\in\Phi_f$. Analogously, we can write $[\Psi]=\{\Psi_f\}_{f\in Aut(\mathcal{M})}$.
\item For each flag $\Phi_f$ we can choose two simple paths $\gamma_f^R$ and $\gamma_f^S$ in $S$ connecting $x_f$ to $x_{fR}$ and $x_f$ to $x_{fS}$ respectively, and such that: $\gamma_f^R\cap\gamma_f^S=x_f$, $\gamma_f^R\cap int(\Psi_f)\neq\emptyset$ and $\gamma_f^S\cap int(\Psi_{fR^{p-1}S})\neq\emptyset$.
\end{enumerate}
In figure \ref{FIGURE4} we illustrate the definition of $i:G\hookrightarrow S$.
As it can be seen, on vertices $i(f)=x_f$ and on edges $i((f,fP))=\gamma_f^R$ and $i((f,fS))=\gamma_f^S$. It is clear that our choice of $i$ is a topological embedding.
%
\begin{remark}
	\label{ob_3}
\begin{upshape}
Let $U$ be a connected open  subset of $G$ with compact boundary and, without loss of generality, such that $\partial U\subset V(G)$. Define:
\[
\widetilde{U}:= Int \left(\bigcup_{f\in U\cap V(G)} (\Phi_f\cup\Psi_{f}\cup\Psi_{fR^{n-1}S})\right)\subset S
\]
\end{upshape}
\end{remark}
and
\begin{equation}\label{ec:6}
\begin{array}{ccccc}
i_{\ast} & : & Ends(G) & \to & Ends(S) \\
          &   & [U_n]_{n\in \N}& \to &[\widetilde{U}_n]_{n\in \N}
\end{array}
\end{equation}
The map $i_\ast$ in the expression \ref{ec:6} is a bijection and the proof is verbatim the same as in the regular case. \qed


\subsection{Proof theorem \ref{TH:endsAutGroups}} Every automorphism group of a regular map with Schl\"afli type $\{p,q\}$ is a quotient of the Coxeter group:
\begin{equation}
	\label{E:CoxeterGroup}
\langle\rho_0,\rho_1,\rho_0\hspace{.5mm}|\hspace{.5mm} \rho_i^2,\hspace{.5mm}(\rho_0\rho_2)^2,\hspace{.5mm}(\rho_0\rho_1)^p,\hspace{.5mm}(\rho_1\rho_2)^q\rangle		 
\end{equation}
and every  automorphism group of a chiral map is a quotient of the Schwarz-group:
\begin{equation}
	\label{E:CoxeterGroup1}
\langle R,S\hspace{.5mm}|\hspace{.5mm} R^p,\hspace{.5mm} S^q,\hspace{.5mm}(RS)^2 \rangle .		
\end{equation}
Each of these groups has Serre's FA property, see~\cite{Serre03} for a definition and a proof of this fact. Serre's FA property is stable when passing to quotients (see \S 6.3.2, [\emph{Ibid}]). Any group having Serre's FA property cannot split as an amalgamated product or an HNN extension. On the other hand, a classical result due to Stallings~\cite{Stallings68} says that if $G$ is finitely generated group with more than one end, then $G$ is either split as an amalgamated product or an HNN extension.\qed



\subsection{Proof of theorems \ref{Coro:RegularMapsMonster} and \ref{Coro:ChiralMapsMonster}}

\label{SUBSEC:ProofCorollaries}

We first address theorem \ref{Coro:RegularMapsMonster}. It is sufficient to show that if the surface $S$ is orientable, $Ends(S)$ is a singleton and $S$ has genus then $S$ has necessarily infinite genus. Given that $S$ has genus, there exists a subsurface $S'\subset S$ with compact closure which is homeomorphic to the punctured torus $S_{1,1}$.  Let $D(S')$ be the saturation of $S'$ by $\mathcal{F}$, that is, $D(S')$ is the (finite) union of all flags in $\mathcal{F}$ that intersect $S'$.
Since both $Ends(S)$ and $Ends(Aut(\mathcal{M}))$ are singletons, $Aut(\mathcal{M})$ is infinite and acts transitively on the set of flags $\mathcal{F}$. Therefore there exists an infinite sequence $\{f_n\}_{n\in\N}$ in $Aut(\mathcal{M})$ such that $f_n(D(S'))\cap f_m(D(S'))=\emptyset$ for all $m\neq n$ in $\N$. This implies that there is an infinite sequence of subsurfaces $S'_n\subset S$ with compact closure, each homeomorphic to $S_{1,1}$ and disjoint by pairs. Hence, $S$ has infinite genus.

We address now theorem \ref{Coro:ChiralMapsMonster}.
There are no chiral maps on the plane. Indeed, let $\mathcal{M}(\Gamma,i,S)$ be a map with Schl\"afli type $\{p,q\}$ on a simply connected orientable surface $S$ with only one end. There are two cases. First when $\frac{1}{p}+\frac{1}{q}=\frac{1}{2}$, case in which one can check per hand that there are no chiral maps for the finitely many possible values of $p$ and $q$. In the second case $\frac{1}{p}+\frac{1}{q}<\frac{1}{2}$. Here $\mathcal{M}(\Gamma,i,S)$ must be a quotient of the universal map in the hyperbolic disk with Schl\"afli type $\{p,q\}$ by some isometry group $H$. Since $S$ is simply connected $H$ must be trivial and then $\mathcal{M}(\Gamma,i,S)$ is  a universal map in the disc with Schl\"afli type $\{p,q\}$. But there are not such chiral maps on the hyperbolic disc, for they always contain reflections.

On the other hand, there are no chiral maps on non-orientable surfaces.  Hence if $\mathcal{M}$ is a chiral map on a surface $S$ such that $Ends(S)$ is a singleton, then $S$ must have genus and $Aut(\mathcal{M})$ must be infinite. As in the regular case let $S'\subset S$ be a subsurface with compact closure homeomorphic to the compact torus $S_{1,1}$ and $D(S')\subset S$ be the saturation of $S'$ by $\mathcal{F}$. This last set must contain elements from both $Aut(\mathcal{M})$-orbits $[\Phi]$ and $[\Psi]$. As in the regular case,
using the transitivity of the action of $Aut(\mathcal{M})$ in one of these two orbits we can extract an infinite sequence $\{f_n\}_{n\in\N}$ in $Aut(\mathcal{M})$ such that $f_n(D(S'))\cap f_m(D(S'))=\emptyset$ for all $m\neq n$ in $\N$. \qed


\end{document}